\documentclass{article}
\usepackage{mathrsfs, amssymb, amsmath}

\begin{document}

\title{Entropies and the derivatives of some Heun functions}
\author{Ioan Ra\c{s}a}

\date{}
\maketitle

\textbf{Abstract}

This short note contains a list of new results concerning the R\'{e}nyi entropy, the Tsallis entropy, and the Heun functions associated with positive linear operators.

\section{Introduction}

This short note contains a list of new results supplementing the articles~\cite{10}, \cite{11}, \cite{12}.

The R\'{e}nyi entropy and the Tsallis entropy associated with positive linear operators have been investigated in~\cite{11} and~\cite{12}. Section~\ref{sect:2} is concerned with two new examples in this direction.

The entropies are naturally related to some Heun functions, as explained in the mentioned articles. Using results from~\cite{4} and~\cite{13}, the derivatives of these Heun functions are studied in Sections~\ref{sect:3} and~\ref{sect:4}.

Detailed proofs will be presented in a forthcoming paper.

\section{Integral operators\label{sect:2}}

Let $B_{m-1}\left ( x_0, x_1, \dots x_m; \cdot \right )$ be the $B$-spline function of degree $m-1$ associated to the equidistant points $x_0<x_1<\dots < x_m$. Consider a given function $\sigma \in \mathscr{C}(\mathbb{R})$ such that $\sigma (x) > 0$, $x \in \mathbb{R}$. For $x, t \in \mathbb{R}$, let
\begin{equation*}
W_n(x,t):=B_{n-1}\left ( x-\sigma (x), x-\sigma (x) + \frac{2\sigma (x)}{n}, \dots , x+\sigma (x);t \right ), \quad n\geq 1.
\end{equation*}

The operator $L_n : \mathscr{C}(\mathbb{R}) \to  \mathscr{C}(\mathbb{R})$ will be defined by
\begin{equation}
L_nf(x):= \int _{\mathbb{R}} W_n (x,t)f(t)dt, \quad n\geq 1, f\in  \mathscr{C}(\mathbb{R}), x\in \mathbb{R}. \label{eq:2.1}
\end{equation}

Several properties of operators of this form are investigated, e.g., in~\cite{1}, \cite{2}, \cite{9}.

The R\'{e}nyi entropy and the Tsallis entropy associated with $L_n$ are, respectively, $-\log{\int _\mathbb{R} W_n^2 (x,t)dt}$ and $1-\int _\mathbb{R} W_n^2 (x,t)dt$; see, e.g., \cite{11} and the references therein. For the operator described by~\eqref{eq:2.1}, we have
\begin{equation*}
\int _\mathbb{R} W^2_n(x,t)dt = \frac{c_n}{\sigma (x)}, \quad x\in \mathbb{R},
\end{equation*}
where $c_n >0$ is a constant depending only on $n$; in particular $c_1 = \frac{1}{2}$, $c_2 = \frac{2}{3}$, $c_3 = \frac{33}{40}$.

Let $e_i(x)=x^i$, $i=0,1,2$. The variance associated with $L_n$ is defined by $V_n (x) = L_n e_2(x) - \left ( L_ne_1(x) \right )^2$, $x\in \mathbb{R}$; see~\cite{11}. In our case,
\begin{equation*}
V_n(x) = \frac{\sigma ^2(x)}{3n}, \quad x \in \mathbb{R}.
\end{equation*}

Remark that the variance, the R\'{e}nyi entropy and the Tsallis entropy are synchronous functions of $x$.

Let now $b_{n,j}(x):= {n \choose j} x^j (1-x)^{n-j}$, $j=0, 1, \dots , n$, $x\in [0,1]$. Then $B_n:\mathscr{C}[0,1] \to \mathscr{C}[0,1]$, $B_nf(x):=\sum _{j=0}^n f \left ( \frac{j}{n} \right ) b_{n,j}(x)$, are the classical Bernstein operators.

The Kantorovich modifications of $B_n$ are defined by
\begin{equation*}
Q_n ^{[k]}f:= \frac{n^k(n-k)!}{n!}D^kB_nf^{(-k)}, \quad f\in \mathscr{C}[0,1], n\geq k \geq 0,
\end{equation*}
where $D$ is the differentiation operator and $f^{(-k)}$ is an antiderivative of order $k$ of $f$; see~\cite{3} and the references therein.

It can be proved that
\begin{equation*}
Q_n^{[k]}f = \sum _{j=0}^{n-k} b_{n-k,j} \int _0^1 f(t)B_{k-1}\left ( \frac{j}{n}, \frac{j+1}{n}, \dots , \frac{j+k}{n} ;t\right ) dt.
\end{equation*}

Therefore, in order to compute the R\'{e}nyi entropy and the Tsallis entropy associated with $Q_n^{[k]}$, we need
\begin{equation*}
S_n^{[k]}(x) := \int _0^1 \left ( \sum _{j=0}^{n-k} b_{n-k,j} (x) B_{k-1}\left ( \frac{j}{n}, \frac{j+1}{n}, \dots , \frac{j+k}{n};t \right ) \right )^2 dt.
\end{equation*}

In fact, it can be proved that
\begin{eqnarray}
S_{n+2}^{[2]}(x) &=& \frac{n+2}{3(n+1)4^n} \sum _{i=0}^n (3n-2i+2)4^i {2i\choose i} {2n-2i\choose {n-i}}\left ( x-\frac{1}{2}\right )^{2i} \nonumber\\
&=& \frac{n+2}{3\pi} \int _0 ^\pi \left ( 1-4x(1-x)\sin ^2{\frac{\phi}{2}} \right )^n \left ( 1+2\cos ^2{\frac{\phi}{2}}\right )d\phi .\label{eq:2.2}
\end{eqnarray}

The variance associated with $Q_{n+2}^{[2]}$ is
\begin{equation*}
V_{n+2}^{[2]}(x) = \frac{n}{(n+2)^2}x(1-x) + \frac{1}{6(n+2)^2} , \quad x\in [0,1].
\end{equation*}

Again it follows that the variance, the R\'{e}nyi entropy and the Tsallis entropy associated with $Q_{n+2}^{[2]}$ are synchronous functions.

\section{Heun functions and their derivatives}\label{sect:3}

With classical notation for Heun functions $Hl(a, q; \alpha, \beta ; \gamma, \delta; x)$ and hypergeometric functions $_2F_1 (a,b;c;x)$ we have (see~\cite{12})

\begin{eqnarray}
Hl \left ( \frac{1}{2}, q; 2q, 1; 1, 1; x \right ) &=& (1-x)^{-2q} {_2F_1} \left ( q,q; 1; \left ( \frac{x}{x-1} \right )^2\right )\nonumber \\ &=& \frac{1}{\pi} \int _0^\pi \left ( 1-4x(1-x)\sin ^2{\frac{\phi}{2}} \right )^{-q}d\phi. \label{eq:3.1}
\end{eqnarray}

The integral representation can be compared with~\eqref{eq:2.2}.

Combined with
\begin{equation*}
_2F_1(a,b;c;x) = (1-x)^{-a}{_2F_1} \left ( a,c-b;c;\frac{x}{x-1}\right ),
\end{equation*}
\eqref{eq:3.1} leads to
\begin{equation}
Hl \left ( \frac{1}{2}, q;2q,1;1,1;x \right ) = (1-2x)^{-q} {_2F_1} \left ( q,1-q;1;\frac{x^2}{2x-1}\right ). \label{eq:3.2}
\end{equation}

With notation from~\cite{10}, \cite{12}, let
\begin{equation*}
F_n(x) : = \sum _{k=0}^n \left ( {n \choose k} x^k (1-x)^{n-k} \right )^2,
\end{equation*}
\begin{equation*}
G_n(x) : = \sum _{k=0}^\infty \left ( {n +k-1\choose k} x^k (1+x)^{-n-k} \right )^2,
\end{equation*}
\begin{equation*}
U_n(x) : = \sum _{k=0}^n \left ( {n \choose k} x^k (1+x)^{-n} \right )^2,
\end{equation*}
\begin{equation*}
J_n(x) : = \sum _{k=0}^\infty \left ( {n+k \choose k} x^k (1-x)^{n+1} \right )^2.
\end{equation*}

It was proved in~\cite{12} that
\begin{equation}
F_n(x) = Hl \left (\frac{1}{2}, -n; -2n, 1; 1, 1; x \right ),\label{eq:3.3}
\end{equation}
\begin{equation}
G_n(x) = Hl \left (\frac{1}{2}, n; 2n, 1; 1, 1; -x \right ),\label{eq:3.4}
\end{equation}
and, moreover, the polynomial Heun function $F_n(x)$ and the rational Heun function $G_n(-x)$ are related by
\begin{equation}
G_n(-x) = (1-2x)^{1-2n}F_{n-1}(x). \label{eq:3.5}
\end{equation}

Similarly, $U_n(x)$ and $J_n(x)$ can be expressed as
\begin{equation}
U_n(x) = F_n \left ( \frac{x}{x+1}\right ), \label{eq:3.6}
\end{equation}
\begin{equation}
J_n(x) = \left ( \frac{1-x}{1+x} \right )^{2n+1} F_n \left ( \frac{1}{1-x} \right ).\label{eq:3.7}
\end{equation}

Recall that the Legendre polynomials $P_n(x)$ are related to $_2F_1$ by
\begin{equation}
_2F_1 (-n,n+1;1;x) = P_n(1-2x).\label{eq:3.8}
\end{equation}

From~\eqref{eq:3.3}, \eqref{eq:3.2} and~\eqref{eq:3.8} we get
\begin{equation}
F_n(x) = (1-2x)^n P_n \left ( \frac{2x^2-2x+1}{1-2x}\right ). \label{eq:3.9}
\end{equation}

This formula was proved (with a different method) by Thorsten Neuschel~\cite{7} and Geno Nikolov~\cite{8}; it was used in order to prove a conjecture involving the polynomials $F_n(x)$.

The derivative of a Heun function $Hl (a,q;\alpha, \beta; \gamma , \delta ;x)$ satisfying
\begin{equation}
q=a\alpha \beta \label{eq:3.10}
\end{equation}
was studied in~\cite{4}. From the corresponding results we infer, for $\gamma \neq 0, -1, -2, \dots ,$
\begin{eqnarray}
&&\frac{d}{dx}Hl \left ( \frac{1}{2}, \frac{1}{2}\alpha \beta; \alpha , \beta ; \gamma , \gamma ; x \right )\label{eq:3.11} \\ & =& \frac{\alpha \beta}{\gamma} (1-2x)Hl \left ( \frac{1}{2}, \frac{1}{2}(\alpha +2)(\beta +2);\alpha +2, \beta +2; \gamma +1 , \gamma +1; x \right ), \nonumber
\end{eqnarray}

\begin{eqnarray}
&&\frac{d}{dx}Hl \left ( \frac{1}{2}, \frac{1}{2}\alpha \beta; \alpha , \beta ; \gamma , \gamma ; x \right )\label{eq:3.12} \\ & =& \frac{\alpha \beta}{\gamma} (1-2x)^{2\gamma -\alpha -\beta -1}Hl \left ( \frac{1}{2}, \frac{1}{2}(2\gamma-\alpha)(2\gamma - \beta );2\gamma - \alpha , 2\gamma -\beta ; \gamma +1 , \gamma +1; x \right ). \nonumber
\end{eqnarray}

The equality between the right-hand sides of~\eqref{eq:3.11} and~\eqref{eq:3.12} follows also from line 3 in Table 2 of~\cite{6}. Let us remark that the Heun functions in these right-hand sides satisfy also the condition~\eqref{eq:3.10}, so that it is possible to express their derivatives in terms of other Heun functions.

From~\eqref{eq:3.3} and~\eqref{eq:3.11} we get
\begin{equation}
\frac{d}{dx}F_n(x) = 2n(2x-1)Hl \left ( \frac{1}{2}, 3-3n;2-2n,3;2,2;x \right ),\label{eq:3.13}
\end{equation}
and finally, for $i=0,1, \dots ,n$,
\begin{equation}
Hl \left ( \frac{1}{2}, (i-n)(2i+1);2(i-n), 2i+1;i+1,i+1;x \right )\label{eq:3.14}
\end{equation}
\begin{equation*}
= \frac{(2i)!!}{(2i-1)!!}4^{-n}{n \choose i}^{-1} \sum _{j=0}^{n-i} 4^j {i+j \choose i} {2i+2j \choose i+j} {2n-2i-2j \choose n-i-j} \left ( x-\frac{1}{2} \right )^{2j}.
\end{equation*}

To conclude this section, let us remark that the Heun function from~\eqref{eq:3.1} satisfies the condition (33) from~\cite{5}; the consequence of this fact will be investigated elsewhere.

\section{A confluent Heun function}\label{sect:4}

Let $u(x) = HC (p, \gamma , \delta, \alpha, \sigma; x)$ be the solution of the confluent Heun equation
\begin{equation}
u''(x) + \left ( 4p + \frac{\gamma}{x} + \frac{\delta}{x-1} \right ) u'(x) + \frac{4p\alpha x - \sigma}{x(x-1)}u(x) = 0, \label{eq:4.1}
\end{equation}
with $u(0)=1$ (See~\cite{13}). From~\cite[(21)]{13} we get
\begin{equation}
\frac{d}{dx}HC(p,\gamma, 0, \alpha, 4p\alpha;x) = -\frac{\sigma}{\gamma} HC (p, \gamma +1, 0, \alpha +1, 4p(\alpha +1);x), \label{eq:4.2}
\end{equation}
\begin{equation}
\frac{d}{dx}HC(p,\gamma, 0, \alpha, 4p\alpha;x) = \frac{\sigma}{\gamma}(x-1)HC(p,\gamma +1, 2, \alpha +2, 4p (\alpha +1)-\gamma -1; x). \label{eq:4.3}
\end{equation}

With notation from~\cite{10}, \cite{12}, let
\begin{equation*}
K_n(x):=\sum _{k=0}^\infty \left ( e^{-nx}\frac{(nx)^k}{k!} \right )^2.
\end{equation*}

Then (see~\cite{10}, \cite{12}),
\begin{equation}
xK_n''(x) + (4nx+1)K_n'(x) + 2nK_n(x) = 0.\label{eq:4.4}
\end{equation}

We get immediately
\begin{equation}
K_n(x) = HC \left ( n,1,0,\frac{1}{2},2n;x \right ).\label{eq:4.5}
\end{equation}

From~\eqref{eq:4.3} and~\eqref{eq:4.5} it follows that
\begin{equation}
HC \left ( n, 2, 2, \frac{5}{2}, 6n-2;x \right ) = \frac{1}{2n(x-1)}K_n'(x).\label{eq:4.6}
\end{equation}

Similarly, from~\eqref{eq:4.2} and~\eqref{eq:4.5},
\begin{equation}
HC \left ( n, 2, 0, \frac{3}{2}, 6n;x \right ) = -\frac{1}{2n}K_n'(x).\label{eq:4.7}
\end{equation}

Now applying repeatedly~\eqref{eq:4.2} we can generalize~\eqref{eq:4.7} to
\begin{equation}
HC \left ( n, j+1, 0, \frac{2j+1}{2}, 2n(2j+1);x \right ) = \frac{K_n^{(j)}(x)}{K_n^{(j)}(0)}, \quad j \geq 0, \label{eq:4.8}
\end{equation}
where
\begin{equation}
K_n^{(j)}(0) = (-2n)^j \sum _{i=0}^{[j/2]} {j \choose 2i} {2i \choose i}4^{-i}.\label{eq:4.9}
\end{equation}

From~\eqref{eq:4.4} we get
\begin{equation}
\left ( K_n^{(j)}\right )'' + \left (4n+\frac{j+1}{x} \right )\left ( K_n^{(j)}\right )' + \frac{2n(2j+1)}{x}K_n^{(j)} = 0,\label{eq:4.10}
\end{equation}
and this provides an alternative proof of~\eqref{eq:4.8}.

\end{document}